\documentclass[11pt]{article}
\usepackage{amsmath, amssymb, amscd}
\usepackage[matrix,arrow]{xy}
\def\eqref#1{(\ref{#1})}
\def\Bbb#1{\mathbb #1}

\newcommand{\goth}{\frak}

\newcommand{\arrow}{{\:\longrightarrow\:}}

\newcommand{\C}{{\Bbb C}}
\newcommand{\R}{{\Bbb R}}

\renewcommand{\H}{{\Bbb H}}
\newcommand{\6}{\partial}
\renewcommand{\1}{\sqrt{-1}\:}

\newcommand{\cntrct}                
{\hspace{2pt}\raisebox{1pt}{\text{$\lrcorner$}}\hspace{2pt}}

\makeatletter
\def\x@arrow{\DOTSB\Relbar}
\def\xlongequalsignfill@{\arrowfill@\x@arrow\Relbar\x@arrow}
\newcommand{\xlongequal}[2][]{%
\ext@arrow 0099\xlongequalsignfill@{#1}{#2}} \makeatother

\def\RR{\mathbb{R}}
\def\CC{\mathbb{C}}

\def\ZZ{\mathbb{Z}}
\def\HH{\mathbb{H}}

\def\pt{\partial}

\renewcommand{\tilde}{\widetilde}
\renewcommand{\bar}{\overline}
\renewcommand{\phi}{\varphi}
\renewcommand{\epsilon}{\varepsilon}
\renewcommand{\geq}{\geqslant}
\renewcommand{\leq}{\leqslant}

\newcommand{\End}{\operatorname{\text{\sf End}}}

\newcommand{\Hol}{\operatorname{Hol}}

\newcommand{\const}{\operatorname{const}}

\newcommand{\comment}[1]{{}}
\def\blacksquare{\hbox{\vrule width 4pt height 4pt depth 0pt}}
\def\endproof{\blacksquare}

\makeatletter \@ifundefined{Bbb}
    {\newcommand{\Bbb}[1]{{\mathbb #1}}}%
{}

\newcommand{\ps@verbit}{%
 \renewcommand{\@oddhead}{%
         \scriptsize
         {Calabi problem for HKT-manifolds}
         \hfil\tiny {S. Alesker, M. Verbitsky, \ \ \ \ \@date}}
 \renewcommand{\@evenhead}{\@oddhead}
 \renewcommand{\@oddfoot}{\hfil\thepage\hfil}
 \renewcommand{\@evenfoot}{\@oddfoot}}
\newcounter{Mycounter}[section]
\newcounter{lemma}[section]
\renewcommand{\thelemma}{{Lemma \thesection.\arabic{lemma}}}
\newcommand{\lemma}{%
    \setcounter{lemma}{\value{Mycounter}}
    \refstepcounter{lemma}
    \stepcounter{Mycounter}
    {\bf \thelemma:\ }}
\newcounter{claim}[section]
\newcounter{sublemma}[section]
\newcounter{corollary}[section]
\newcounter{theorem}[section]
\renewcommand{\thetheorem}{{Theorem \thesection.\arabic{theorem}}}
\newcommand{\theorem}{%
    \setcounter{theorem}{\value{Mycounter}}
    \refstepcounter{theorem}
    \stepcounter{Mycounter}
    {\bf \thetheorem:\ }}
\newcounter{conjecture}[section]
\renewcommand{\theconjecture}{{Conjecture \thesection.\arabic{conjecture}}}
\newcommand{\conjecture}{%
    \setcounter{conjecture}{\value{Mycounter}}
    \refstepcounter{conjecture}
    \stepcounter{Mycounter}
    {\bf \theconjecture:\ }}
\newcounter{proposition}[section]
\renewcommand{\theproposition}
      {{Proposition \thesection.\arabic{proposition}}}
\newcommand{\proposition}{%
    \setcounter{proposition}{\value{Mycounter}}
    \refstepcounter{proposition}
    \stepcounter{Mycounter}
    {\bf \theproposition:\ }}
\newcounter{definition}[section]
\renewcommand{\thedefinition}
      {{Definition~\thesection.\arabic{definition}}}
\newcommand{\definition}{%
    \setcounter{definition}{\value{Mycounter}}
    \refstepcounter{definition}
    \stepcounter{Mycounter}
    {\bf \thedefinition:\ }}

\newcounter{example}[section]
\newcounter{remark}[section]
\renewcommand{\theremark}{{Remark \thesection.\arabic{remark}}}
\newcommand{\remark}{%
    \setcounter{remark}{\value{Mycounter}}
    \refstepcounter{remark}
    \stepcounter{Mycounter}
    {\bf \theremark:\ }}

\newcounter{problem}[section]
\newcounter{question}[section]
\@addtoreset{equation}{section} \@addtoreset{footnote}{section}
\makeatother

\begin{document}

\begin{center}
{\LARGE\bf Quaternionic Monge-Amp\`ere equation\\[2mm] and
Calabi problem for HKT-manifolds }

\hfill

S. Alesker \footnote{ Semyon Alesker was partially supported by ISF
grant 1369/04.}, M. Verbitsky,\footnote{Misha Verbitsky is an EPSRC
advanced fellow supported by CRDF grant RM1-2354-MO02 and EPSRC
grant  GR/R77773/01}
\end{center}

\def\ome{\omega}
\def\Ome{\Omega}
{\small \hspace{0.15\linewidth}

\begin{minipage}[t]{0.8\linewidth}
{\bf Abstract}\\ A quaternionic version of the Calabi problem on the
Monge-Amp\`ere equation is introduced, namely a quaternionic
Monge-Amp\`ere equation on a compact hypercomplex manifold with an
HKT-metric. The equation is non-linear elliptic of second order. For
a hypercomplex manifold with holonomy in $SL(n,{\Bbb H})$,
uniqueness (up to a constant) of a solution is proven, as well as
the zero order a priori estimate. The existence of a solution is
conjectured, similar to the Calabi-Yau theorem. We reformulate this
quaternionic equation as a special case of the complex Hessian
equation, making sense on any complex manifold.
\end{minipage}

}
{ \small \tableofcontents }
\section{Introduction}\label{S:introduction}

We introduce a quaternionic version of the Calabi problem on the
Monge-\-Amp\`ere equation. The problem is motivated by close analogy
(discussed below) with the classical Calabi-Yau theorem
\cite{_Yau:Calabi-Yau_} on complex Monge-Amp\`ere equations on
K\"ahler manifolds. Our version of the Calabi problem is about the
quaternionic Monge-Amp\`ere equation on a compact hypercomplex
manifold with an HKT-metric. The equation is non-linear elliptic of
second order. When a holonomy group of the Obata connection of a
hypercomplex manifold lies in $SL(n, {\Bbb H})$, uniqueness (up to a
constant) of the solution is proven, as well as the zero order a
priori estimate. Under this assumption, we conjecture existence of
the solution. We give a reformulation of this quaternionic equation
as a special case of the complex Hessian equation, which makes sense
on any complex manifold.

\hfill

\begin{definition}
A {\itshape hypercomplex} manifold is a smooth manifold $M$ together
with a triple $(I,J,K)$ of complex structures satisfying the usual
quaternionic relations:
$$IJ=-JI=K.$$
\end{definition}

Necessarily the real dimension of a hypercomplex manifold is
divisible by 4. The simplest example of a hypercomplex manifold is
the flat space $\HH^n$.

\hfill

\begin{remark}
(1) In this article we will presume (unlike in much of the
literature on the subject) that the complex structures $I,J,K$ act
on the {\itshape right} on the tangent bundle $TM$ of $M$. This
action extends uniquely to the right action of the algebra $\HH$ of
quaternions on $TM$.

(2) It follows that the dimension of a hypercomplex manifold $M$ is
divisible by 4.

(3) Hypercomplex manifolds were introduced explicitly by Boyer
\cite{boyer}.
\end{remark}

\hfill

Let $M$ be a hypercomplex manifold, and $g$ a Riemannian metric on
$M$. The metric $g$ is called {\itshape quaternionic Hermitian} (or
hyperhermitian) if $g$ is invariant with respect to the group
$SU(2)\subset \HH^*$ of unitary quaternions. Given a quaternionic
Hermitian metric $g$ on a hypercomplex manifold $M$, consider the
differential form
$$\Ome:=-\ome_J+\sqrt{-1}\ome_K$$
where $\ome_L(A,B):=g(A,B\circ L)$ for any $L\in \HH$ with $L^2=-1$
and any vector fields $A,B$ on $M$. It is easy to see that $\Ome$ is
a $(2,0)$-form with respect to the complex structure $I$.

\begin{definition}\label{def-hkt-metr}
The metric $g$ on $M$ is called {\bf an HKT-metric} (here HKT stands
for HyperK\"ahler with Torsion) if
$$\6 \Ome =0$$
where $\6$ is the usual $\6$-differential on the complex manifold
$(M,I)$
\end{definition}

The form $\Omega$ corresponding to an HKT-metric $g$ will be called
{\bf an HKT-form.}

\begin{remark}
HKT-metrics on hypercomplex manifolds first were introduced by Howe
and Papadopoulos \cite{howe-papa}. Their original definition was
different but equivalent to \ref{def-hkt-metr} (see
\cite{_Gra_Poon_}).
\end{remark}

HKT-metrics on hypercomplex manifolds are analogous in many respects
to K\"ahler metrics on complex manifolds. For example, it is a
classical result that any K\"ahler form $\omega$ on a complex
manifold can be locally written in the form $\omega=dd^ch$ where $h$
is a strictly plurisubharmonic function called a potential of
$\omega$, and vice versa (see, e.g. \cite{_Griffi_Harri_}).
Similarly, by \cite{_Banos_Swann_} (see also
\cite{_Alesker_Verbitsky_}), an HKT-form $\Omega$ on a hypercomplex
manifold locally admits a potential: it can be written as
$$\Omega=\6\6_JH$$ where $\6_J=J^{-1}\circ \bar \6 \circ J$,
and $H$ is a strictly plurisubharmonic function in the quaternionic
sense. The converse is also true. The notion of quaternionic
plurisubharmonicity is relatively new: on the flat space $\HH^n$ it
was introduced in \cite{alesker-bsm-03} and independently by G.
Henkin around the same time (unpublished), and on general
hypercomplex manifolds in \cite{_Alesker_Verbitsky_}. More recently
the notion of plurisubharmonic functions has been generalized to yet
another context of calibrated geometries
\cite{harvey-lawson-psh-calib}.

Motivated by the analogy with the complex case, we introduce the
following quaternionic version of the Calabi problem. Let
$(M^{4n},I,J,K)$ be a compact hypercomplex manifold of real
dimension $4n$. Let $\Omega$ be an HKT-form. Let $f$ be a
real-valued $C^\infty$ function on $M$. The quaternionic Calabi
problem is to study solvability of the following quaternionic
Monge-Amp\`ere equation with an unknown real-valued function $\phi$:
\begin{eqnarray}\label{mampere}
(\Omega+\6\6_J\phi)^n=e^f\Omega^n.
\end{eqnarray}

By \ref{L:interior} below, if  a $C^\infty$-function $\phi$
satisfies the Monge-Amp\`ere equation (\ref{mampere}) then
$\Omega+\6\6_J\phi$ is an HKT-form, namely it corresponds to a new
HKT-metric. This equation is a non-linear elliptic equation of
second order. We formulate the following conjecture.

\hfill

\begin{conjecture}\label{Conj:CY}
Let us assume that $(M,I)$ admits a holomorphic (with respect to the
complex structure $I$) non-vanishing $(2n,0)$-form $\Theta$. Then
the quaternionic Monge-Amp\`ere equation (\ref{mampere}) has a
$C^\infty$-solution $\phi$ provided the following necessary
condition on the initial data is satisfied:
$$\int_M(e^f-1)\Omega^n\wedge\bar\Theta=0.$$
\end{conjecture}

\hfill

In this article we show that under the condition of existence of
such $\Theta$ a solution of (\ref{mampere}) is unique up to a
constant (\ref{C:ellip-unique}). Our next main result is a zero
order a priori estimate (\ref{ze-6}): there exists a constant $C$
depending on $M,\Omega$, and $||f||_{C^0}$ only, such that the
solution $\phi$ satisfying the normalization condition
$\int_M\phi\cdot \Omega^n\wedge \bar \Theta=0$ must satisfy the
estimate
$$||\phi||_{C^0}\leq C,$$
where $||\cdot||_{C^0}$ denotes the maximum norm on $M$, i.e.
$||u||_{C^0}:=\max\{|u(x)|\, \,|\, x\in M\}$. Our proof of this
estimate is a modification of Yau's argument \cite{_Yau:Calabi-Yau_}
in the complex case as presented in \cite{joyce}.

\hfill

\begin{remark}
Let us comment on how restrictive the condition of existence of a
form $\Theta$ is. Recall that a hypercomplex manifold $M$ carries a
unique torsion free connection such that the complex structures
$I,J,K$ are parallel with respect to it. It it called the Obata
connection as it was discovered by Obata \cite{obata}. It was shown
by the second named author \cite{_Verbitsky:Canonical_HKT_} that if
$M$ is a compact HKT-manifold admitting a holomorphic (with respect
to $I$) $(2n,0)$-form $\Theta$, then the holonomy of the Obata
connection is contained in the group $SL_n(\HH)$ (instead of
$GL_n(\HH)$). Conversely, if the holonomy of the Obata connection is
contained in $SL_n(\HH)$, then there exists a form $\Theta$ as above
which, moreover, can be chosen to be q-positive (in sense of Section
\ref{Ss:R_V} below).
\end{remark}

\begin{remark}
The quaternionic Monge-Amp\`ere equation (\ref{mampere}) can be
interpreted in the following geometric way. Assume we are given an
HKT-form $\Ome$, and a strongly q-real $(2n,0)$-form on $(M,I)$ (see
Section \ref{Ss:R_V} for the definition) which is nowhere vanishing
and hence may be assumed to have the form $e^f\Ome^n$. We are
looking for a new HKT-form of the form $\Ome+\6\6_J\phi$ whose
volume form is equal to the prescribed form $e^f\Ome^n$.
\end{remark}

\hfill

We note that the Calabi problem also has its real version where it
becomes a real Monge-Amp\`ere equation on smooth compact manifolds
with an affine flat structure. This real Calabi problem was first
considered and successfully solved by Cheng and Yau
\cite{cheng-yau}. Note also that the classical Dirichlet problem for
the Monge-Amp\`ere equation in strictly pseudoconvex domains has its
quaternionic version considered and partly solved by the first named
author \cite{alesker-jga-03}. We refer to \cite{alesker-jga-03} for
the details.

Finally, in this article we present a reformulation of the
quaternionic Monge-Amp\`ere equation as a special case of a complex
Hessian equation on the complex manifold $X$, $\dim_\CC X=m$. This
equation is:
\begin{equation}\label{general_mampere}
(\omega-\sqrt{-1}\6\bar\6 \phi)^n\wedge\Phi=e^f\omega^n\wedge\Phi,
\end{equation}
where $n=m-k$, $\Phi\in \Lambda^{k,k}(X)$, $\omega \in
\Lambda^{1,1}(X)$, $f\in C^\infty(X)$ are fixed, and $\Phi, \omega$
satisfy some positivity assumptions (see \ref{P:first _form}). Let
us state the conditions more explicitly under an additional
assumption of existence of the form $\Theta$ as in \ref{Conj:CY}.

We consider the Monge-Amp\`ere equation (\ref{general_mampere})
where the unknown function $\phi$ belongs to the class of $C^\infty$
functions, such that $\omega-\sqrt{-1}\6\bar\6 \phi$ lies in the
interior of the cone of $\Phi$-positive forms.
\ref{_CY_MA_quat_Theorem_} states that the quaternionic
Monge-Amp\`ere equation (\ref{mampere}) is equivalent to
(\ref{general_mampere}) for appropriate choices of $\Phi$ and
$\omega$ under the assumption of existence of the form $\Theta$ as
in \ref{Conj:CY}. Moreover  one may assume
$d\Phi=d(\Phi\wedge\omega)=0$ (see
\ref{_P:Phi-character_Proposition_}). We show that if all these
conditions are satisfied on a complex compact manifold $X$, then the
complex Hessian equation (\ref{general_mampere}) is elliptic, its
solution is unique up to a constant, and a necessary condition for
solvability is
$$\int_X(e^f-1)\omega^n\wedge \Phi=0.$$
We refer to \ref{C:8} below for the details.


\section{Quaternionic Dolbeault complex}
\label{_q.D._Section_}


To continue, we need a definition and some properties of the Salamon
complex on hypercomplex manifolds. The following Section is adapted
from \cite{_Verbitsky:qD_}. The quaternionic cohomology is a
well-known subject, introduced by S. Salamon
(\cite{_Capria-Salamon_}, \cite{_Salamon_}, \cite{_Baston_},
\cite{_Leung_}). Here we give an exposition of quaternionic
cohomology and a quaternionic Dolbeault complex for hypercomplex
manifolds.


\subsection{Quaternionic Dolbeault complex: the definition}

Let $M^{4n}$ be a hypercomplex manifold of real dimension $4n$, and
\[ \Lambda^0 (M) \stackrel d\arrow \Lambda^1 (M)
  \stackrel d\arrow \Lambda^2 (M) \stackrel d\arrow ...
\]
its de Rham complex. Consider the natural (left) action of $SU(2)$
on $\Lambda^*M$. Clearly, $SU(2)$ acts on $\Lambda^i(M)$, $i\leq
\frac 1 2 \dim_\R M$ with weights \[ i, i-2, i-4, ... \] We denote
by $\Lambda^i_+(M)$ the maximal $SU(2)$-subspace of $\Lambda^i(M)$,
on which $SU(2)$ acts (on the left) with weight $i$. We again
emphasize that necessarily $i\leq 2n= \frac{1}{2}\dim_\RR M$.

\hfill

The following linear algebraic lemma allows one to compute
$\Lambda^i_+(M)$ explicitly.

\hfill

\lemma\label{_Lambda_+_explicit_Lemma_} (\cite[Proposition
2.9]{_Verbitsky:qD_}) With the above assumptions, let $I$ be the
induced complex structure, and $\H_I$ the quaternion space,
considered as a 2-dimensional complex vector space with the complex
structure induced by $I$ when $I$ acts on $\HH_\CC$ on the right.
Denote by $\Lambda^{p,0}_I(M)$ the space of the $(p,0)$-form on $M$,
with respect to the Hodge decomposition associated with the complex
structure $I$. The space $\H_I$ is equipped with the natural action
of $SU(2)$. Consider $\Lambda^{p,0}_I(M)$ as a representation of
$SU(2)$, with trivial group action. Then, there is a canonical
isomorphism
\begin{equation}\label{_Lambda_+_explicit_Equation_}
\Lambda^p_+(M,\CC) \cong S^p_\C \H_I \otimes_\C \Lambda^{p,0}_I(M),
\end{equation}
where $S^p_\C \H_I$ denotes a $p$-th symmetric power of $\H_I$ over
$\CC$. Moreover, the $SU(2)$-action on $\Lambda^p_+(M)$ is
compatible with the isomorphism \eqref{_Lambda_+_explicit_Equation_}
when $SU(2)$ acts trivially on $\Lambda^{p,0}_I(M)$.

\hfill

{\bf Proof:} Fix a standard basis $1, I, J, K$ in $\H$. Since $\HH$
acts on the tangent bundle $TM$ on the right, $\HH$ acts on
$\Lambda^1(M)$ on the left, namely we have a canonical map
$$\HH\otimes _\RR\Lambda^1(M)\to \Lambda^1(M).$$
Taking the complexification, we get a $\CC$-linear map
$$\HH\otimes _\RR\Lambda^1(M,\CC)\to \Lambda^1(M,\CC).$$
Since $\Lambda^{1,0}(M)\subset\Lambda^1(M,\CC)$, we get the map
\begin{eqnarray}\label{map1}
\HH\otimes _\RR\Lambda^{1,0}(M)\to \Lambda^1(M,\CC).
\end{eqnarray}

We have a canonical quotient map
$$\HH\otimes _\RR\Lambda^{1,0}(M)\to \HH_I\otimes_\CC\Lambda^{1,0}(M).$$
It is easy to see that the map (\ref{map1}) factorizes uniquely via
a map
\begin{eqnarray}\label{map2}
\HH_I\otimes_\CC\Lambda^{1,0}(M)\to \Lambda^1(M,\CC).
\end{eqnarray}

(It is easy to write down the map (\ref{map2}) explicitly. Let $h_1,
h_2\in \H_I$ be the basis in $\H_I$: $h_1=1$, $h_2= J$. Then
$h_1\otimes x \mapsto x,\, h_2\otimes x\mapsto J(x)$.)
 Consider the $SU(2)$-equivariant homomorphism
\begin{equation}\label{_basis_in_H_I_Equation_}
 \H_I \otimes_\CC \Lambda^{1,0}_I(M)\arrow \Lambda^1(M),
\end{equation}
mapping $h_1\otimes \eta$ to $\eta$ and $h_2\otimes \eta$ to
$J(\eta)$, where $J$ denotes an endomorphism of $\Lambda^1(M)$
induced by $J$. The isomorphism \eqref{_Lambda_+_explicit_Equation_}
is obvious for $p=1$:
\begin{equation}\label{_Lambda_+^1_Equation_}
\Lambda^1(M)=\Lambda^1_+(M)= \H_I \otimes_\CC \Lambda^{1,0}_I(M)
\end{equation}
This isomorphism is by construction $SU(2)$-equivariant. Given two
vector spaces $A$ and $B$, we have a natural map
\begin{equation}\label{_exte_symme_Equation_}
 S^iA\otimes \Lambda^iB\arrow \Lambda^i(A\otimes B),
\end{equation}
given by
$$(a_1\otimes \dots \otimes a_i)\otimes (b_1\wedge \dots\wedge b_i)\mapsto
\frac{1}{i!}\sum_{\sigma\in \Sigma_i} (a_{\sigma(1)}\otimes
b_1)\wedge \dots \wedge (a_{\sigma(i)}\otimes b_i).$$
 {}From \eqref{_exte_symme_Equation_} and
\eqref{_Lambda_+^1_Equation_}, we obtain the natural
$SU(2)$-equivariant map
\begin{equation}\label{_Lambda_+_explicit_homomo_interme_Equation_}
S^p_\C \H_I \otimes_\C \Lambda^{p,0}_I(M) \arrow \Lambda^p(M,\CC).
\end{equation}
Since $S^p_\C \H_I$ has weight $p$, the arrow
\eqref{_Lambda_+_explicit_homomo_interme_Equation_} maps $S^p_\C
\H_I \otimes_\C \Lambda^{p,0}_I(M)$ to $\Lambda^p_+(M)$. We have
constructed the map
\begin{equation}\label{_Lambda_+_explicit_homomo_Equation_}
S^p_\C \H_I \otimes_\C \Lambda^{p,0}_I(M)
  \stackrel \Psi\arrow \Lambda^p_+(M).
\end{equation}
It remains to show that it is an isomorphism. Let
$adI:\Lambda^*M\arrow \Lambda^*M$ act on the $(p,q)$-forms
$ad(\eta)=(p-q)\1\eta$. Clearly, $-\1ad I$ is a root of the Lie
algebra $SU(2)$. It is well known that an irreducible representation
of a Lie algebra is generated by a highest weight vector. For the
Lie algebra $\goth{su}(2)$, this means that $\Lambda^p_+(M)$ is a
subspace of $\Lambda^p(M)$ generated by $SU(2)$ from the subspace
$W\subset \Lambda^p_+(M)$ consisting of all vectors on which $-\1ad
I$ acts as a multiplication by $p$. On the other hand, $W$ coincides
with $\Lambda^{p,0}_I(M)$. We obtained the following:
\begin{equation}\label{_p,0_genera_Lambda_+_Equation_}
\begin{minipage}[m]{0.8\linewidth}
The space $\Lambda^p_+(M)$ is generated by $SU(2)$ from its subspace
$\Lambda^{p,0}_I(M)$.
\end{minipage}
\end{equation}
The image of
\[
\Psi:\; S^p_\C \H_I \otimes_\C \Lambda^{p,0}_I(M) \arrow
\Lambda^p_+(M)
\]
 is an $SU(2)$-invariant subspace of $\Lambda^p(M)$
containing $\Lambda^{p,0}_I(M)$. By
\eqref{_p,0_genera_Lambda_+_Equation_}, this means that $\Psi$ is
surjective. Let $R\subset S^p_\C \H_I \otimes_\C \Lambda^{p,0}_I(M)$
be the kernel of $\Psi$. By construction, $R$ is $SU(2)$-invariant,
of weight $p$. By the same arguments as above, $R$ is generated by
its subspace of highest weight, that is, the vectors of type
$h_1^p\eta$, where $\eta \in \Lambda^{p,0}_I(M)$ (see
\eqref{_basis_in_H_I_Equation_}). On the other hand, on the subspace
\[ h_1^p\cdot \Lambda^{p,0}_I(M)\subset
 S^p_\C \H_I \otimes_\C \Lambda^{p,0}_I(M),
\]
the map $\Psi$ is, by construction, injective. Therefore, the
intersection $h_1^p\cdot \Lambda^{p,0}_I(M)\cap R$ is zero. We have
proved that $\Psi$ is an isomorphism.
\ref{_Lambda_+_explicit_Lemma_} is proven.
\endproof

\hfill

Consider an $SU(2)$-invariant decomposition
\begin{equation}\label{_decompo_Lambda_to_Lambda_+_Equation_}
\Lambda^p(M) = \Lambda^p_+(M)\oplus V^p
\end{equation}
where $V^p$ is the sum of all $SU(2)$-subspaces of $\Lambda^p(M)$ of
weight less than $p$. Since $SU(2)$-action is multiplicative on
$\Lambda^*(M)$, the subspace $\tilde V:=\oplus_p V^p\subset
\Lambda^*(M)$ is an ideal. Therefore, the quotient
\[
  \Lambda^*_+(M)= \Lambda^*(M)/\tilde V
\]
is an algebra. Using the decomposition
\eqref{_decompo_Lambda_to_Lambda_+_Equation_}, we define the
quaternionic Dolbeault differential $d_+:\; \Lambda^*_+(M)\arrow
\Lambda^*_+(M)$ as the composition of the de Rham differential and
the projection $\Lambda^*(M)\to \Lambda^*_+(M)$. Since de Rham
differential cannot increase the $SU(2)$-weight of a form by more
than 1, $d$ preserves the subspace $V^*\subset \Lambda^*(M)$.
Therefore, $d_+$ is a differential in $\Lambda^*_+(M)$.

\hfill

\definition
Let
\[ \Lambda^0 (M) \stackrel {d_+}\arrow \Lambda^1 (M)
  \stackrel {d_+}\arrow \Lambda^2_+ (M) \stackrel {d_+}\arrow
  \Lambda^3_+ (M) \stackrel {d_+}\arrow...\arrow \Lambda^{2n}_+(M)
\]
be the differential graded algebra constructed above\footnote{We
identify $\Lambda^0 M$ and $\Lambda^0_+ M$, $\Lambda^1 M$ and
$\Lambda^1_+ M$.}. It is called {\bf the quaternionic Dolbeault
complex}, or {\bf Salamon complex}.

\hfill

\remark The isomorphism \eqref{_Lambda_+_explicit_Equation_} is
clearly multiplicative:
$$\oplus_p S_\C^p(\H_I)\otimes \Lambda^{p,0}_I(M)\simeq \oplus_p
\Lambda^p_+(M,\C).$$
Notice that, in the course of the proof of
\ref{_Lambda_+_explicit_Lemma_}, we have proven the following result
(see \ref{_p,0_genera_Lambda_+_Equation_}).

\hfill

\begin{claim}\label{Cl:inclusion}
For any $p\geq 0$
$$\Lambda^{p,0}_I(M)\subset \Lambda^p_+(M).$$
\end{claim}


\subsection[Hodge decomposition for the quaternionic Dolbeault complex]
{Hodge decomposition for the quaternionic  \\ Dolbeault complex}
\label{_Hodge_Dolbeault_Subsection_}

Let $M$ be a hypercomplex manifold and $I$ an induced complex
structure.  As usually, we have the operator $adI:\; \Lambda^*(M)
\arrow \Lambda^*(M)$ mapping a $(p,q)$-form $\eta$ to $\1(p-q)\eta$.
By definition, $ad I$ belongs to the Lie algebra $\goth{su}(2)$
acting on $\Lambda^*(M)$ in the standard way. Therefore, $ad I$
preserves the subspace $\Lambda^*_+(M) \subset \Lambda^*(M)$. We
obtain the Hodge decomposition
\[ \Lambda^*_+(M) = \oplus_{p+q\leq 2n} \Lambda^{p,q}_{+,I}(M). \]

\hfill

\definition\label{_Hodge_for_q.D._Definition_}
The decomposition
\[ \Lambda^*_+(M) = \oplus_{p+q\leq 2n} \Lambda^{p,q}_{+,I}(M)
\]
is called {\bf the Hodge decomposition for the quaternionic
Dolbeault complex}.

\hfill

\subsection{The Dolbeault bicomplex and quaternionic Dolbeault
complex} \label{_qD-bico_Subsection_}
Let $M^{4n}$ be a hypercomplex manifold of real dimension $4n$ and
$I, J, K\in \H$ the standard triple of induced complex structures.
Clearly, $J$ acts on the complexified cotangent space
$\Lambda^1(M,\C)$ mapping $\Lambda_I^{0,1}(M)$ to
$\Lambda_I^{1,0}(M)$. Consider a differential operator
\[ \6_J:\; C^\infty(M)\arrow \Lambda_I^{1,0}(M),\]
mapping $f$ to $J^{-1}(\bar\6 f)$, where $\bar\6:\;
C^\infty(M)\arrow \Lambda_I^{0,1}(M)$ is the standard Dolbeault
differential on the complex manifold $(M,I)$. We extend $\6_J$ to a
differential
\[
  \6_J:\; \Lambda_I^{p,0}(M)\arrow \Lambda_I^{p+1,0}(M),
\]
using the Leibnitz rule. Then $\pt_J=J^{-1}\circ \bar\pt\circ J$.

\hfill


%

\hfill

\proposition \label{_d_+_Hodge_6_J_Proposition_} (see also
\cite[Theorem 2.10]{_Verbitsky:qD_}) Let $M^{4n}$ be a hypercomplex
manifold, $I$ an induced complex structure, $I, J, K$ the standard
basis in quaternion algebra, and
\[ \Lambda^*_+(M)= \oplus_{p+q\leq 2n} \Lambda^{p,q}_{I,+}(M) \]
the Hodge decomposition of the quaternionic Dolbeault complex
(Subsection \ref{_Hodge_Dolbeault_Subsection_}). Then there exists a
canonical isomorphism
\begin{equation}\label{_Lambda_+_Hodge_deco_expli_Equation_}
  \Lambda^{p,q}_{I, +}(M)\cong \Lambda^{p+q, 0}_I(M).
\end{equation}
Under this identification, the quaternionic Dolbeault differential
\[ d_+:\; \Lambda^{p,q}_{I, +}(M)\arrow
  \Lambda^{p+1,q}_{I, +}(M)\oplus \Lambda^{p,q+1}_{I, +}(M)
\]
corresponds to the sum
\[ \6 \oplus \6_J:\; \Lambda^{p+q, 0}_{I}(M) \arrow
  \Lambda^{p+q+1, 0}_{I}(M)\oplus \Lambda^{p+q+1,0}_{I}(M).
\]
{\bf Proof:} Consider the isomorphisms
\eqref{_Lambda_+_explicit_Equation_}
\begin{equation}\label{_Lambda_+_expli_to_Hodge_Equation_}
\Lambda^p_+(M,\C) \cong S^p_\C \H_I\otimes _\C \Lambda^{p,0}_I(M).
\end{equation}
The Hodge decomposition of
\eqref{_Lambda_+_expli_to_Hodge_Equation_} is induced by the
$SU(2)$-action, as follows. Let $\rho_I:\  U(1)\arrow SU(2)$ be the
group homomorphism defined by
$\rho_I(e^{\sqrt{-1}\theta})=e^{I\theta}$ for any $\theta\in
\RR/2\pi\ZZ$. From the definition of the $SU(2)$-action, it follows
that the Hodge decomposition of $\Lambda^*(M)$ coincides with the
weight decomposition under the action of $\rho_I:\; U(1) \arrow
\End(\Lambda^*(M))$. The $SU(2)$-action on $S^p_\C \H_I\otimes _\C
\Lambda^{p,0}_I(M)$ is trivial on the second component. Consider the
weight decomposition
\[ S^i_\C \H_I\cong \bigoplus\limits_{p+q=i}S^{p,q}_\C \H_I
\]
associated with $\rho_I$. Then
\eqref{_Lambda_+_expli_to_Hodge_Equation_} translates to the
isomorphism
\[ \Lambda^{p,q}_{I, +}(M) \cong S^{p,q}_\C \H_I\otimes_\C
\Lambda^{p+q,0}_I(M).
\]
Let $h_1$, $h_2$ be the basis in $\H_I$ defined as in the proof of
\ref{_Lambda_+_explicit_Lemma_}, i.e. $h_1=1,\, h_2=J$. An
elementary calculation shows that $h_1$ has weight (1,0), and $h_2$
has weight (0,1). Therefore, the space $S^{p,q}_\C \H_I$ is
1-dimensional and generated by $h_1^ph_2^q$. We have obtained an
isomorphism
\begin{equation}\label{_h_1_2_q-D_decompo_Equation_}
 \Lambda^{p,q}_{I, +}(M) \cong h_1^p\cdot h_2^q\cdot \Lambda^{p+q,0}_I(M).
\end{equation}
This proves \eqref{_Lambda_+_Hodge_deco_expli_Equation_}. The
isomorphism  \eqref{_h_1_2_q-D_decompo_Equation_} is multiplicative
by Remark \eqref{_Lambda_+_explicit_Equation_}. Consider the
differential
\[ \hat d_+ = h_1 \6 +h_2 \6_J :\;
  S^p_\C \H_I\otimes _\C \Lambda^{p,0}_I(M) \arrow
  S^{p+1}_\C \H_I\otimes _\C \Lambda^{p+1,0}_I(M)
\]
To prove our proposition, we need to show that the quaternionic
Dolbeault differential $d_+$ coincides with $\hat d_+$ under the
identification \eqref{_h_1_2_q-D_decompo_Equation_}. The isomorphism
\eqref{_h_1_2_q-D_decompo_Equation_} is multiplicative, and the
differentials $d_+$ and $\hat d_+$ both satisfy the Leibnitz
rule.\footnote{The differential $\hat d_+$ satisfies the Leibnitz
rule, because $\pt\pt_J=-\pt_J\pt$. The last equation  is clear: the
differentials $d$, $d_I:= -I d I, d_J:= -JdJ, d_K:= -KdK$
anticommute because of integrability of $I,J,K$.} Therefore, it is
sufficient to show that
\begin{equation}\label{_diffe_equal_Equation_}
d_+ = \hat d_+
\end{equation}
on $C^\infty(M) = \Lambda^0_+(M)$. On functions, the equality
\eqref{_diffe_equal_Equation_} is immediately implied by the
definition of the isomorphism
\[ \Lambda^{1}_+(M)\cong \H_I\otimes _\C \Lambda^{1,0}_I(M).
\]
\ref{_d_+_Hodge_6_J_Proposition_} is proven.\endproof

\hfill

The statement of \ref{_d_+_Hodge_6_J_Proposition_} can be
represented by the following diagram:
\begin{equation}\label{_bicomple_XY_Equation}
\begin{minipage}[m]{0.85\linewidth}
{\tiny $ \xymatrix @C+1mm @R+10mm@!0  {
 && \Lambda^0_+(M) \ar[dl]^{d'_+} \ar[dr]^{d''_+}
  && && && \Lambda^{0,0}_I(M) \ar[dl]^{h_1 \6} \ar[dr]^{h_2 \6_J}
  &&  \\
 & \Lambda^{1,0}_+(M) \ar[dl]^{d'_+} \ar[dr]^{d''_+} &
 & \Lambda^{0,1}_+(M) \ar[dl]^{d'_+} \ar[dr]^{d''_+}&&
\text{\large $\cong$} &
 &h_1\Lambda^{1,0}_I(M)\ar[dl]^{h_1 \6} \ar[dr]^{h_2 \6_J}&  &
 h_2\Lambda^{1,0}_I(M)\ar[dl]^{h_1 \6} \ar[dr]^{h_2 \6_J}&\\
 \Lambda^{2,0}_+(M) && \Lambda^{1,1}_+(M)
  && \Lambda^{0,2}_+(M)& \ \ \ \ \ \ & h_1^2\Lambda^{2,0}_I(M)& &
h_1h_2\Lambda^{2,0}_I(M) & & h_2^2\Lambda^{2,0}_I(M) \\
} $ }
\end{minipage}
\end{equation}
where $d_+= d'_+ + d''_+$ is the Hodge decomposition of the
quaternionic Dolbeault differential.

\hfill

\definition
With the above assumptions, the bicomplex
\eqref{_bicomple_XY_Equation} is called {\bf the quaternionic
Dolbeault bicomplex}.

\hfill


\hfill

\begin{lemma}\label{L:projection2}
The projection of $\eta\in\Lambda^{1,1}_I(M)$ to the
$SU(2)$-invariant part of $\Lambda^2(M)$ is given by the map
\begin{eqnarray}\label{eq:01}
\eta \arrow \frac{1}{2}\left(\eta(\cdot, \cdot)+\eta(\cdot \circ
J,\cdot \circ J)\right).
\end{eqnarray}

\end{lemma}

{\bf Proof:} It is easy to see that the 2-form (\ref{eq:01}) is
invariant under $I$ and $J$, hence under $K$. This implies that this
2-form is $SU(2)$-invariant. Also if $\eta$ was already
$SU(2)$-invariant then (\ref{eq:01}) is equal to $\eta$. \endproof

\hfill

\begin{lemma}\label{L:omega}
Let $g$ be a quaternionic Hermitian metric on a hypercomplex
manifold $(M,I,J,K)$. Define
$$\omega_I(X,Y):=g(X,Y\circ I).$$
Then $\omega_I\in \Lambda^{1,1}_+(M)$.
\end{lemma}

{\bf Proof:} It is clear that $\omega_I\in \Lambda^{1,1}_I(M)$.
Since \[ \Lambda^2(M)=\Lambda^2_+(M)\oplus \Lambda^2_{SU(2)}(M),\]
to prove the lemma we have to check that the projection of
$\omega_I$ to the $SU(2)$-invariant forms vanishes. By
\ref{L:projection2} this projection is equal to
$$\frac{1}{2}\left(g(X,Y\circ I)+g(X\circ J,Y\circ JI)\right)=0.$$
The lemma is proven.
\endproof

\hfill

\begin{lemma}\label{L:Omega_omega}
Let $g$ be a quaternionic Hermitian metric on a hypercomplex
manifold $(M,I,J,K)$. Define $\omega_I\in \Lambda^{1,1}_+(M)$ as in
\ref{L:omega} and
$$\Omega(X,Y):=-(g(X,Y\circ J)-\sqrt{-1}g(X,Y\circ K)).$$
Then under the isomorphism \eqref{_Lambda_+_explicit_Equation_}
$(h_1\cdot h_2)\otimes\Omega$ corresponds to $\sqrt{-1}\omega_I$.
\end{lemma}

{\bf Proof:} Under the isomorphism
\ref{_Lambda_+_explicit_Equation_} the form $(h_1\otimes
h_2)\otimes\Omega$ corresponds to the form $\zeta$ given by
\begin{align*}
\zeta(X,Y)&=\frac{1}{2}(\Omega(X,Y\circ J)+\Omega(X\circ J,Y))=\\
\ &-\frac{1}{2}\bigg(g(X,Y\circ J^2)-\sqrt{-1}g(X,Y\circ JK)\\
\ & +g(X\circ J,Y\circ J)-\sqrt{-1}g(X\circ J,Y\circ K)\bigg)=\\
\ & =\sqrt{-1}g(X,Y\circ I)=\sqrt{-1}\omega_I(X,Y).
\end{align*}

The lemma is proven. \endproof
\section{Quaternionic Monge-Amp\`ere equation}


\subsection[First reformulation of quaternionic Monge-Amp\`ere equation.]{First reformulation of quaternionic Monge-\\Amp\`ere equation.}

\label{_comple_Hess_Equation_}

Let $(M,I,J,K)$ be a hypercomplex manifold. Let $g$ be a
quaternionic Hermitian Riemannian metric. Define as in Section
\ref{_qD-bico_Subsection_} the 2-forms

\begin{align*}
\omega_I(X,Y):=&g(X,Y\circ I)\in \Lambda^{1,1}_+(M),\\
\Omega(X,Y):=&-(g(X,Y\circ J)-\sqrt{-1}g(X,Y\circ K)).
\end{align*}

We want to rewrite the quaternionic Monge-Amp\`ere equation
\begin{equation}\label{eq:11}
(\Omega+ \6 \6_J \phi)^n = A e^f \Omega^n
\end{equation}
in terms of the $\Lambda_+^{*,*}(M)$-bicomplex.
Let us multiply both sides of \eqref{eq:11} by $h_1^{n}\cdot
h_2^{n}$ and apply the isomorphism
\eqref{_Lambda_+_explicit_Equation_}. We get
\begin{equation} \label{_q_M-a_in_terms_of_qD_Equation_}
P_+((\sqrt{-1}\omega_I+ d'_+d''_+ \phi)^n) = P_+(A e^f
(\sqrt{-1}\omega_I)^n)
\end{equation}
where $P_+:\; \Lambda^*(M) \arrow \Lambda^*_+(M)$ is a natural
$SU(2)$-invariant projection. Indeed the isomorphism
\eqref{_Lambda_+_explicit_Equation_} is multiplicative and, by
\ref{L:Omega_omega}, carries $h_1\cdot h_2\cdot\Omega$ to
$\sqrt{-1}\omega_I$ and $\6\6_J$ to $d'_+d''_+$. But it is easy to
see that $d'_+d''_+=P_+\circ \6\bar\6$. Hence equation \eqref{eq:11}
is equivalent to the equation
\begin{eqnarray}\label{eq:12}
P_+((\omega_I-\sqrt{-1}\6\bar\6\phi)^n)=Ae^fP_+(\omega_I^n).
\end{eqnarray}

\hfill

\begin{lemma}\label{L:lemma on product}
For any $\eta\in \Lambda^{2n}_{+}(M)$, any non-negative integer $m$,
and any $\xi\in \Lambda^m(M)$,
$$\eta\wedge\xi=\eta\wedge P_+(\xi).$$
\end{lemma}

{\bf Proof:} It is enough to show that if $\xi\in \Lambda^m(M)$
belongs to a subspace of $SU(2)$-weight at most $m-1$ then
$\eta\wedge\xi=0$. In this case, the Clebsch-Gordan formula implies
that $\eta\wedge\xi$ belongs to a subspace of $\Lambda^{2n+m}(M)$
generated by $SU(2)$-weights $2n+m-1,2n+m-3,\dots,2n-m+1$. But
$\Lambda^{2n+m}(M)$ has no vectors for these weights because it is
dual to $\Lambda^{2n-m}(M)$, all of the $SU(2)$-weights of the
latter space are less than or equal to $2n-m$. \endproof

\hfill


\begin{proposition}\label{P:first _form}
Let $M$ be a hypercomplex quaternionic Hermitian manifold, and
\begin{equation}\label{_H-C-Y_standard_Equation_}
(\Omega+ \6 \6_J \phi)^n = A e^f \Omega^n,
\end{equation}
the quaternionic Monge-Amp\`ere equation. Then
\eqref{_H-C-Y_standard_Equation_} is equivalent to the following
equation
\begin{equation}\label{_H-C-Y_via_P_Equation_}
 (\omega_I-\sqrt{-1} \6\bar\6 \phi)^n\wedge P_+(\omega_I^n) =
 A e^f \omega_I^n\wedge P_+(\omega_I^n).
\end{equation}
\end{proposition}

{\bf Proof:} We need to check that \eqref{_H-C-Y_via_P_Equation_} is
equivalent to \eqref{eq:12}. Both sides of \eqref{eq:12} belong to
$\Lambda^{n,n}_+(M)$. However, the vector space
$\Lambda^{n,n}_+(M)\cong \Lambda^{2n,0}(M)$ is clearly
1-dimensional, and generated by $P_+(\omega_I^n)$. By \ref{L:lemma
on product}, for any $\eta,\xi\in \Lambda^{2n}(M)$, one has
$P(\eta)\wedge\xi=\eta\wedge P(\xi)=P(\eta)\wedge P(\xi)$. This
implies that the equation \eqref{eq:12} is equivalent to this
equation multiplied by $\omega_I^n$, and the latter is equivalent to
\begin{eqnarray}\label{eq:13}
(\omega_I-\sqrt{-1}\6\bar\6\phi)^n\wedge
P_+(\omega^n)=Ae^f\omega_I^n\wedge P_+(\omega_I^n).
\end{eqnarray}

 \endproof

\hfill

\remark Equation \eqref{_H-C-Y_via_P_Equation_} is a special case of
the so-called {\bf complex Hessian equation}. More generally, a
generalized complex Hessian equation is written as
$\Psi(\sqrt{-1}\6\bar\6 u) =f$, where $\Psi$ is a symmetric
polynomial in the eigenvalues of the $(1,1)$-form $\sqrt{-1}\6\bar\6
u$.

\hfill


\subsection{Operators $R$ and $V$}\label{Ss:R_V}
We would like to present yet another reformulation of the
quaternionic Monge-Amp\`ere equation. For this we introduce, in this
section, two operators $R$ and $V$ on differential forms. Denote by
\[ \tilde R:\; \Lambda^{p,q}_{I,+}(M)=h_1^ph_2^q\Lambda^{p+q,0}_I(M)
\tilde\arrow \Lambda^{p+q,0}_{I}(M)
\]
the isomorphism constructed in \ref{_d_+_Hodge_6_J_Proposition_}.
Let \[ R:\; \Lambda^{p,q}_{I}(M) \arrow \Lambda^{p+q,0}_{I}(M)\] be
the composition of the standard projection
\[ \Lambda^{p,q}_{I}(M)\overset{P_+}{\arrow} \Lambda^{p,q}_{I,+}(M)\]
with $\tilde R$. In \cite{_Alesker_Verbitsky_}, we defined a real
structure on $\Lambda^{2p,0}_{I}(M)$, that is, an anticomplex
involution mapping $\lambda\in \Lambda^{2p,0}_{I}(M)$ into
$J\bar\lambda$ (since $I$ and $J$ anticommute, $J$ maps
$(p,q)$-forms into $(q,p)$-forms). Forms fixed under this involution
we call {\bf q-real} (q stands for quaternions). We also define a
notion of positivity: a real $(2,0)$-form $\eta$ is {\bf q-positive}
if $\eta (X,X\circ J)\geq 0$ for any real vector field $X$. A {\bf
strongly q-positive cone} is the cone of q-real $(2p,0)$-forms which
is generated by the products of positive forms with non-negative
coefficients (this definition is parallel to one given in complex
analysis - see e.g. \cite{_Demailly_}). It can be shown that this
convex cone is closed and has non-empty interior. A q-real
$(2p,0)$-form $\eta$ is called {\bf weakly q-positive} if for any
strongly q-positive $(2n-2p,0)$-form $\xi$ the product
$\eta\wedge\xi\in \Lambda^{2n,0}_I(M)$ is strongly q-positive. The
set of weakly q-positive forms is a closed convex cone with
non-empty interior. Note that any strongly q-positive form is weakly
q-positive, and the notions of weak and strong q-positivity coincide
for $(0,0)$, $(2n,0)$, $(2,0)$, and $(2n-2,0)$-forms (see
\cite{alesker-adv-05}, Propositions 2.2.2 and 2.2.4, where the
q-positivity is called just positivity, and only the flat space
$M=\HH^n$ is considered). The map $R$ satisfies the following
properties.

\hfill

\theorem\label{_R_properties_Theorem_} (see also \cite[Claim 4.2,
Claim 4.5]{_Verbitsky:positive_}) Let $M$ be a hypercomplex manifold
and
\[ R:\; \oplus_{p,q}\Lambda^{p,q}_{I}(M) \arrow
\oplus_{p,q}(h_1^ph_2^q)\otimes\Lambda^{p+q,0}_{I}(M)=\oplus_rS_\C^r(\H_I)\otimes
\Lambda^{r,0}_I(M)
\] the map constructed above. Then

\begin{description}
\item[(i)] $R$ is multiplicative: $R(x\wedge y) = R(x)\wedge R(y)$.
\item[(ii)] $R$ is related to the real structures as follows:
\[
R(\bar\lambda) =(-1)^p \bar{J(R(\lambda))},
\]
for any $\lambda\in \Lambda^{p,q}_I(M)$ where the action of $J$ on
$S^r_\CC(\HH_I)$ is identical.
\item[(iii)]
We have
\[
R(\6\lambda) = \6 R(\lambda),\ \  R(\bar\6\lambda) = \6_J
R(\lambda).
\]
\item[(iv)] $(\sqrt{-1})^pR$ maps strongly positive $(p,p)$-forms (in the
complex sense) to strongly positive $(2p,0)$-forms (in the
quaternionic sense).
\end{description}

{\bf Proof:} \ref{_R_properties_Theorem_} (i) is clear from the
construction.

\def\ome{\omega}
\def\Ome{\Omega}
\def\Lam{\Lambda}
Let us prove part (ii). Due to the multiplicativity of $R$ it is
enough to check the statement for $\lambda\in \Lambda^1(M)$.  Set
$\eta:=R(\lambda)\in\Lambda^{1,0}_I(M),\,
\xi:=R(\bar\lambda)\in\Lambda^{1,0}_I(M)$. First assume that
$\lambda\in \Lambda^{1,0}_I(M)$. Then
\begin{eqnarray}\label{E:forms1}
\lambda=h_1\eta=\eta,\ \ \ \bar\lambda=h_2\xi=J(\xi).
\end{eqnarray}
We have to show that $\xi=-\overline{J(\eta)}$ which is obvious by
(\ref{E:forms1}). Let us assume now that $\lambda\in \Lambda^{0,1}$.
We have
\begin{eqnarray}\label{E:forms2}
\lambda=h_2\eta=J(\eta),\ \ \  \bar\lambda=h_1\xi=\xi.
\end{eqnarray}
We have to show that $\xi=\overline{J(\eta)}$ which is obvious by
(\ref{E:forms2}).

Let us prove part (iii). We have
\begin{align*}
R(\6\lambda)=&\tilde R(P_+(\6\lambda))=\tilde
R(P_+(\6(P_+\lambda)))\\
=&\tilde R(d'_+(P_+\lambda))
\xlongequal{\text{\tiny\ref{_d_+_Hodge_6_J_Proposition_}}} \6(\tilde
R(P_+\lambda))=\6(R\lambda).
\end{align*}

Similarly one proves the equality $R(\bar\6\lambda) = \6_J
R(\lambda)$. Let us prove part (iv). Again due to the
multiplicativity of $R$ it is enough to prove it for 2-forms, i.e.
$p=1$. First recall that
$\Lambda^{1,1}_I(M)=\Lambda^{1,1}_{I,+}(M)\oplus
\Lambda^2_{SU(2)}(M)$ where $\Lambda^2_{SU(2)}(M)$ denotes the space
of $SU(2)$-invariant 2-forms (which are necessarily of type (1,1) on
$(M,I)$). Let $\ome\in \Lambda^{1,1}_I(M)$. By \ref{L:projection2}
its projection $P_{SU(2)}(\ome)$ to $\Lambda^2_{SU(2)}(M)$ is equal
to
$$P_{SU(2)}(\ome)(X,Y)=\frac{1}{2}(\ome(X,Y)+\ome(XJ,YJ)).$$
Then the projection $P_+(\ome)$ to $\Lambda^{1,1}_{I,+}(M)$ is equal
to
$$P_+(\ome)(X,Y)=\frac{1}{2}(\ome(X,Y)-\ome(XJ,YJ)).$$
Hence
\begin{align*}
P_+(\ome)(X,XI)&=\frac{1}{2}(\ome(X,XI)-\ome(XJ,XIJ))\\=&
\frac{1}{2}(\ome(X,XI)+\ome(XJ,(XJ)I)).
\end{align*}

It follows that if $\ome$ is positive then $P_+(\ome)$ is positive.
Next we have the equality
$\Lambda^{1,1}_{I,+}(M)=h_1h_2\Lambda^{2,0}_{I}(M)$. It remains to
show that $\Ome\in\Lambda^{2,0}_I(M)$ is positive provided
$\sqrt{-1}h_1h_2\Ome\in \Lambda^{1,1}_{I,+}(M)$ is positive. We have
$$((h_1h_2)\cdot\Omega)(X,X\circ I)=\frac{1}{2}(\Omega(X,X\circ
IJ)+\Omega(X\circ J,X\circ I))=\Omega(X,X\circ K).$$ But for any
$\Omega\in \Lambda^{2,0}_I(M)$ and any vector field $X$ one has
$\Omega(X,X\circ K)=-\sqrt{-1}\Omega(X,X\circ J)$. Hence
$$\sqrt{-1}((h_1h_2)\cdot\Omega)(X,X\circ I)=\Omega(X,X\circ
J).$$ Part (iv) is proven.
\endproof

\hfill

We will need also the following lemma.

\hfill

\begin{lemma}\label{L:su2-positive}
Let $\eta\in \Lambda^{1,1}_I(M)$ be positive (in the complex sense).
If $\eta$ is $SU(2)$-invariant then $\eta=0$.
\end{lemma}

{\bf Proof:} For any real vector field $X$ we have $\eta(X,X\circ
I)\geq 0$. Due to the $J$-invariance of $\eta$ we have
$$0\leq \eta(X,X\circ I)=\eta(X\circ J,(X\circ I)\circ
J)=-\eta(X\circ J,(X\circ J)\circ I)\leq 0.$$ Hence $\eta(X,X\circ
I)=0$ for any real vector field $X$. But since $\eta$ has type
$(1,1)$ this implies that $\eta=0$. \endproof

\hfill

Fix a non-vanishing holomorphic section $\Theta\in
\Lambda^{2n,0}_I(M)$ of the canonical class. Assume moreover that
$\Theta$ is q-real and q-positive. We define a map
\[
V:\; \Lambda^{2p,0}_{I}(M)\arrow \Lambda^{n+p,n+p}_{I}(M)
\]
by the following relation
\begin{equation}\label{_V_via_test_forms_Equation_}
V(\eta) \wedge \xi=\eta \wedge R(\xi) \wedge \bar\Theta,
\end{equation}
where $\xi\in \Lambda^{n-p,n-p}_{I}(M)$ is an arbitrary test form,
and $\eta \in \Lambda^{2p,0}_{I}(M)$.

\hfill

\theorem\label{_V_prop_Theorem_} Let $(M,I,J,K)$ be a hypercomplex
manifold equipped with a non-vanishing holomorphic section
$\Theta\in \Lambda^{2n,0}_I(M)$ of the canonical class. Assume that
$\Theta$ is q-real and q-positive. Then \[ V:\;
\Lambda^{2p,0}_{I}(M)\arrow \Lambda^{n+p,n+p}_{I}(M)\] satisfies the
following properties:

\begin{description}
\item[(i)] For any $\eta\in \Lambda^{2p,0}_I(M)$, one has
$$V(\overline{J\eta})=\overline{V(\eta)}.$$
In particular $V$ maps q-real $(2p,0)$-forms to real (in the complex
sense) $(n+p,n+p)$-forms.
\item[(ii)]
A form $\eta\in \Lambda^{2p,0}_{I}(M)$ is $\6$-exact ($\6$-closed,
$\6_J$-exact, $\6_J$-closed) if and only if $V(\eta)$ is
$\bar\6$-exact ($\bar\6$-closed, $\bar\6_J$-exact, $\bar\6_J$-closed
respectively).

\item[(iii)] $V$ maps weakly q-positive forms  to
weakly positive (in the complex sense) forms.
\item[(iv)] $V:\;
\Lambda^{2p,0}_{I}(M)\arrow \Lambda^{n+p,n+p}_{I}(M)$ is injective.
\end{description}

{\bf Proof:} \ref{_V_prop_Theorem_} follows from
\ref{_R_properties_Theorem_}, by duality. To see that $V$ maps
q-real forms to real forms, we use
\[
 V(J\bar\eta) \wedge \xi=J\bar \eta \wedge R(\xi) \wedge
 \bar\Theta= \bar\eta \wedge J ( R(\xi) \wedge
 \bar\Theta). \]
(The last equation is true, because $J$ acts on volume
 forms trivially.) Since $\Theta$ is q-real, the last expression is equal
to
\begin{align*}
\bar\eta\wedge J(R(\xi))\wedge\Theta
\xlongequal{\text{\ref{_R_properties_Theorem_}(ii)}}
&\bar\eta\wedge\overline{R(\bar\xi)}\wedge\Theta=\overline{\eta\wedge
R(\bar\xi)\wedge\bar\Theta}=\\
=\overline{V(\eta)\wedge\bar\xi}=&\overline{V(\eta)}\wedge\xi.
\end{align*}
Thus we have shown that $V(J\bar\eta) \wedge
\xi=\overline{V(\eta)}\wedge\xi$ for any $\xi$. This proves
\ref{_V_prop_Theorem_} (i). To check positivity of $V(\eta)$, we use
\ref{_R_properties_Theorem_} (iv) (strongly positive forms are dual
to weakly positive). To show that $V$ maps $\6$-closed forms to
$\6$-closed ones, we use
\[
 \int_M V(\6\eta) \wedge \xi=\int_M \6\eta \wedge R(\xi)
 \wedge \bar\Theta =- \int_M \eta \wedge \6 R(\xi)
 \wedge \bar\Theta = -\int_M \eta \wedge  R(\6\xi)
 \wedge \bar\Theta
\]
(the last equation follows from \ref{_R_properties_Theorem_} (iii)).
Then, for any $\6$-closed $\xi$, $\int_M V(\6\eta) \wedge \xi=0$,
hence $V(\6\eta)$ is exact. The converse is also true, because $R$
is injective (\ref{_d_+_Hodge_6_J_Proposition_}). In a similar way
one deduces the rest of statements of (ii) from
\ref{_R_properties_Theorem_} (iii) and injectivity of $R$. Let us
prove (iv). Assume that $\phi\in \Lambda^{2p,0}_I(M)$ belongs to the
kernel of $V$. Then for any $\xi\in \Lambda^{n-p,n-p}_I(M)$ we have
\[0=V(\phi)\wedge\xi=\phi\wedge R(\xi)\wedge\bar \Theta.
\]
But since $R\colon \Lambda^{n-p,n-p}_I(M)\to
\Lambda^{2(n-p),0}_I(M)$ is onto, and $\Theta\in
\Lambda^{2n,0}_I(M)$ is non-vanishing this implies that $\phi=0$.
\endproof

\hfill

The following trivial lemma is used later on in this paper.

\hfill

\lemma\label{_V_multi_Lemma_} In assumptions of
\ref{_V_prop_Theorem_}, the following formula is true
\[
V(R(\eta\wedge \nu))= V(R(\eta))\wedge \nu,
\]
for all $\eta\in \Lambda^{p,p}(M), \nu \in \Lambda^{q,q}(M)$.

\hfill

{\bf Proof:} Since $R$ is multiplicative, we have
\[
V(R(\eta\wedge\nu)) \wedge \xi= R(\eta\wedge\nu) \wedge R(\xi)
\wedge \bar\Theta = R(\eta) \wedge R(\nu\wedge\xi) \wedge \bar\Theta
= V(R(\eta))\wedge\nu \wedge \xi
\]
proving \ref{_V_multi_Lemma_}. \endproof

\hfill

Let us define now
\begin{eqnarray}\label{D:Phi}
\Phi:=V(1)\in\Lambda^{n,n}_I(M).
\end{eqnarray}
The following proposition summarizes the main properties of $\Phi$.

\hfill

\proposition\label{_P:Phi-character_Proposition_} The form $\Phi$
satisfies the following properties:

\begin{description}
\item[(i)] $\Phi\in \Lambda^{n,n}_{I,+}(M)$;
\item[(ii)] $\Phi$ is real in the complex sense, i.e. $\bar\Phi=\Phi$;
\item[(iii)]$\Phi$ is weakly positive.
\item[(iv)] $d\Phi=0$.
\item[(v)] For any Hermitian form $\omega\in
 \Lambda^{1,1}_{I,+}(M)$, the product $\Phi
\wedge \omega^{n-1}$ belongs to the interior of the cone of strongly
(= weakly) positive $(2n-1,2n-1)$-forms.
\item[(vi)]  A Hermitian form $\omega\in \Lambda^{1,1}_{I,+}(M)$ is HKT
if and only if $\Phi \wedge \omega$ is closed. In this case
$\Phi\wedge\omega^j$ is closed for any $j$.
\end{description}

{\bf Proof:} To prove (i) it is enough to show that for any $\xi\in
\Lambda^{n,n}_I(M)$ which belongs to the subspace of $(n,n)$-forms
generated by $SU(2)$-weights at most $2n-1$, one has
$\Phi\wedge\xi=0$. But $\Phi\wedge\xi=R(\xi)\wedge\bar \Theta$, and
$R(\xi)=0$. Thus (i) is proven.

Part (ii) follows immediately from \ref{_V_prop_Theorem_} (i). Part
(iii) follows from \ref{_V_prop_Theorem_} (iii). Let us prove (iv).
\ref{_P:Phi-character_Proposition_} (iv) is clear  from
\ref{_V_prop_Theorem_} (ii), because $\Phi=V(1)$, and 1 is closed.

Let us prove (v). To prove that $\Phi \wedge \omega^{n-1}$ lies in
the interior of the cone of positive elements, let us suppose to the
contrary that it lies on the boundary. Since the cone of (strongly)
positive $(1,1)$-forms is closed there exists $\eta\in
\Lam^{1,1}_I(M)$ such that $\eta\geq 0$, $\eta\ne 0$, and
$$\omega_I^{n-1}\wedge \Phi\wedge\eta=0.$$
But by (\ref{_V_via_test_forms_Equation_})
$$\omega_I^{n-1}\wedge
\Phi\wedge\eta=R(\omega_I^{n-1}\wedge\eta)\wedge\bar\Theta=
(R(\omega_I))^{n-1}\wedge R(\eta)\wedge\bar\Theta=0.$$ Hence
$(R(\omega_I))^{n-1}\wedge R(\eta)=0$. Set
$\Omega:=\sqrt{-1}R(\omega_I)$ be the corresponding HKT
$(2,0)$-form. Then $\Omega$ belongs to the interior of the cone of
strongly positive $(2,0)$-forms in the quaternionic sense. This fact
together with the equality $\Omega^{n-1}\wedge R(\eta)=0$ and the
inequality $\sqrt{-1}R(\eta)\geq 0$ (the latter holds by
\ref{_R_properties_Theorem_}), imply that $R(\eta)=0$. But this
means that $\eta$ is an $SU(2)$-invariant 2-form on $M$. But, since
$\eta\geq 0$, \ref{L:su2-positive} implies that $\eta=0$. This
contradiction finishes the proof of (v). Let us prove (vi). Recall
that $\omega$ is HKT if and only if $R(\omega)$ is $\6$-closed. By
by \ref{_V_prop_Theorem_} (ii), this is equivalent to $\6
V(R(\omega))=0$. However, $V(R(\omega))=\omega \wedge V(1)$, by
\ref{_V_multi_Lemma_}, and $\omega \wedge V(1)$ is a real
$(n+1,n+1)$-form, hence $\omega$ is HKT if and only $\omega \wedge
V(1)$ is closed. Then $\omega^k \wedge V(1)= V(R(\omega)^k)$ is also
closed, because $R(\omega)^k$ is a power of an HKT-form $\Omega$,
and
\[
 \6\Omega^k= k\Omega^{k-1} \wedge \6\Omega =0,
\]
by the Leibnitz identity.
\endproof

\hfill

Now we are ready to give yet another reformulation of the
quaternionic Monge-Amp\`ere equation in complex terms under the
additional assumption that we are given a non-vanishing holomorphic
q-real q-positive form $\Theta\in \Lambda^{2n,0}_I(M)$. Let us fix
an HKT-metric on $M$. Let $\Omega\in \Lambda^{2,0}_I(M)$ and
$\omega_I\in \Lambda^{1,1}_{I,+}(M)$ be the corresponding forms.
Namely
$$\omega_I(X,Y)=g(X,Y\circ I),\, \Omega=\sqrt{-1}R(\omega_I).$$
As previously we denote $\Phi:=V(1)\in \Lambda^{n,n}_{I,+}$. Then we
have

\hfill

\theorem \label{_CY_MA_quat_Theorem_} Let $(M^{4n},I,J,K,g)$ be an
HKT-manifold of real dimension $4n$. Consider the quaternionic
Monge-Amp\`ere equation
\begin{equation}\label{_Theta_q_MA_Theorem_}
 (\Omega+ \6 \6_J \phi)^n = e^f \Omega^n.
\end{equation}
Then \eqref{_Theta_q_MA_Theorem_} is equivalent to the following
equation
\begin{equation}\label{_Theta_q_MA_Theorem_2}
(\omega_I-\sqrt{-1} \6 \bar\6 \phi)^n\wedge \Phi = e^f\omega_I^n
\wedge \Phi.
\end{equation}
{\bf Proof:} It is easy to see that
\[ (\Omega+ \6 \6_J \phi)^n =  e^f \Omega^n
\]
is equivalent to
\[
(\Omega+ \6 \6_J \phi)^n\wedge \bar\Theta =
 e^f \Omega^n\wedge \bar\Theta.
\]
However, $R(\sqrt{-1}\omega+\6\bar\6\phi)= \Omega+ \6 \6_J \phi$ as
follows from \ref{_R_properties_Theorem_}. Therefore
\[ (\Omega+ \6 \6_J \phi)^n\wedge \bar\Theta=
  R(\sqrt{-1}\omega+\6\bar\6\phi)^n \wedge
  \bar\Theta=(\sqrt{-1})^n(\omega-\sqrt{-1}\6\bar\6\phi)^n\wedge\Phi
\]
by definition of $\Phi$. On the other hand
$$\Omega^n\wedge\bar\Theta=R((\sqrt{-1}\omega)^n)\wedge\bar\Theta=(\sqrt{-1})^n\omega^n\wedge\Phi.$$
The result follows.
\endproof



\def\eps{\varepsilon}
\section{Complex Hessian equation.}


The goal of this section is to propose a generalization of the
quaternionic Monge-Amp\`ere equation written in the form
\eqref{_Theta_q_MA_Theorem_2} for any complex manifold $X$. Then,
under appropriate assumptions, satisfied in the HKT-case, we prove
ellipticity of the equation and uniqueness of the solution. The main
results of the section are \ref{C:8} and \ref{C:ellip-unique}.
Throughout this section, we fix a complex manifold $X$ of complex
dimension $m$.

\hfill

\begin{definition}
Let $\Phi\in \Lam^{k,k}(X)$. A form $\eta\in \Lam^{p,p}(X)$ is
called $\Phi$-positive if, for any $\nu\in \Lam^{q,q}(X)$ such that
$\Phi\wedge\nu$ is weakly positive, the form
$\Phi\wedge\nu\wedge\eta$ is weakly positive.
\end{definition}

\hfill

\begin{lemma}\label{L:2}
If $\Phi$ is weakly positive, and $\kappa\in \Lam^{p,p}(X)$ is
strongly positive then $\kappa$ is $\Phi$-positive.
\end{lemma}

{\bf Proof} is obvious. \endproof

\hfill

\begin{lemma}\label{L:3}
(i) The set of $\Phi$-positive $(p,p)$-forms is a convex cone. (ii)
If $\Phi$ is weakly positive then the cone of $\Phi$-positive
$(p,p)$-forms has a non-empty interior.
\end{lemma}

{\bf Proof:} Part (i) is obvious. Part (ii) follows from \ref{L:2}
because the sub-cone of strongly positive forms has a non-empty
interior. \endproof

\hfill

\begin{lemma}\label{L:3.5}
Let $\Phi\in \Lam^{k,k}(X)$ be a weakly positive form. Assume that
$\omega_1,\dots,\omega_r$ are $\Phi$-positive. Then
$\omega_1\wedge\dots \omega_r\wedge\Phi$ is also weakly positive.
\end{lemma}

{\bf Proof:} Since $\Phi=\Phi\wedge1$ is weakly positive and
$\omega_1$ is weakly positive then
$\Phi\wedge1\wedge\omega_1=\Phi\wedge\omega_1$ is weakly positive.
Then continue by induction. \endproof

\hfill

\begin{lemma}\label{L:4}
Let $X$ be a compact complex manifold of complex dimension $m$. Let
$\Phi\in \Lambda^{k,k}(X)$ be a weakly positive form. Let $f\in
C^\infty(X)$, $\omega\in \Lambda^{1,1}(X)$ be real. Denote $n:=m-k$.
Consider the Monge-Amp\`ere equation
\begin{eqnarray}\label{E:ma}
(\omega-\sqrt{-1}\6\bar\6 \phi)^n\wedge\Phi=e^f\omega^n\wedge\Phi.
\end{eqnarray}

(i) Assume that the form $(\omega-\sqrt{-1}\6\bar\6
\phi)^{n-1}\wedge\Phi\in \Lam^{m-1,m-1}(X)$ belongs to the interior
of the cone of (strongly=weakly) positive $(m-1,m-1)$-forms. Then
the Monge-Amp\`ere equation (\ref{E:ma}) is elliptic at $\phi$. (ii)
The Monge-Amp\`ere equation (\ref{E:ma}) has at most unique (up to a
constant) solution in the class of $C^\infty$ functions $\phi$
satisfying the following two conditions:

$\bullet$ $(\omega-\sqrt{-1}\6\bar\6 \phi)^{n-1}\wedge\Phi\in
\Lam^{m-1,m-1}(X)$ belongs to the interior of the cone of positive
$(m-1,m-1)$-forms (strongly or weakly they are the same);

$\bullet$ $\omega-\sqrt{-1}\6\bar\6 \phi$ is $\Phi$-positive.
\end{lemma}

\hfill

{\bf Proof:} (i) The linearization of the equation is
$$\psi\mapsto \psi\wedge(\omega-\sqrt{-1}\6\bar\6
\phi)^{n-1}\wedge\Phi,$$ where $\psi\in \Lam^{1,1}(X)$. This
operator is obviously elliptic. (ii) Let $\phi_1,\phi_2$ be two
solutions as in (ii). Then they satisfy

{\small
\begin{eqnarray*}
dd^c(\phi_1-\phi_2)\wedge\left(\sum_{k=0}^{n-1}
(\omega-\sqrt{-1}\6\bar\6\phi_1)^k\wedge(\omega-\sqrt{-1}\6\bar\6\phi_2)^{n-1-k}\wedge\Phi\right)=0
\end{eqnarray*}
}

By \ref{L:3.5}, the form
$(\omega-\sqrt{-1}\6\bar\6\phi_1)^k\wedge(\omega-\sqrt{-1}\6\bar\6\phi_2)^{n-1-k}\wedge\Phi$
is weakly positive for each $k$. Moreover, for $k=0$, this form
belongs to the interior of the cone of (strongly=weakly) positive
$(m-1,m-1)$-forms. Then the function $\phi_1-\phi_2$ satisfies the
linear elliptic equation of second order on the compact manifold
$X$. Hence it must be constant by the strong maximum principle
(\cite{_Gilbarg_Trudinger_}).
\endproof

\hfill

\begin{lemma}\label{L:7}
Let $X$ be a complex manifold of complex dimension $m$. Let $\Phi\in
\Lam^{k,k}(X)$ be weakly positive. Denote as previously $n=m-k$.
Assume moreover that there exists a (strongly) positive form $\gamma
\in \Lam^{1,1}(X)$ such that $\gamma^{n-1}\wedge\Phi\in
\Lam^{m-1,m-1}(X)$ belongs to the interior of the cone of positive
$(m-1,m-1)$-forms (weakly or strongly they are the same). Let
$\eta\in\Lam^{1,1}(X)$ belong to the interior of the cone of
$\Phi$-positive forms. Then $\eta^{n-1}\wedge\Phi$ belongs to the
interior of the cone of positive $(m-1,m-1)$-forms.
\end{lemma}

\hfill

{\bf Proof:} Multiplying $\gamma$ by a small $\eps>0$, we may assume
that $\eta-\gamma$ is $\Phi$-positive. We have
$$\eta^{n-1}\wedge\Phi=(\gamma+(\eta-\gamma))^{n-1}\wedge\Phi=
\gamma^{n-1}\wedge\Phi+\sum_{j=0}^{n-2}\gamma^j\wedge(\eta-\gamma)^{n-1-j}\wedge\Phi.$$
Every summand in the second sum is (weakly) positive by \ref{L:3.5},
while $\gamma^{n-1}\wedge\Phi$ belongs to the interior of positive
$(m-1,m-1)$-forms. Hence the whole sum also belongs to the interior
of positive $(m-1,m-1)$-forms. \endproof

\hfill

As a corollary we deduce the main result of this section.

\hfill

\begin{theorem}\label{C:8}
Let $X$ be a compact complex manifold of complex dimension $m$. Let
$\Phi\in \Lam^{k,k}(X)$ be a weakly positive form such that there
exists a (strongly) positive form $\gamma \in \Lam^{1,1}(X)$ with
the property that, for $n:=m-k$, the form $\gamma^{n-1}\wedge\Phi$
belongs to the interior of the cone of positive $(m-1,m-1)$-forms
(weakly or strongly they are the same). Let $f\in C^\infty(X)$,
$\omega\in \Lambda^{1,1}(X)$ be real. Consider the Monge-Amp\`ere
equation
\begin{eqnarray}\label{E:ma1}
(\omega-\sqrt{-1}\6\bar\6 \phi)^n\wedge\Phi=e^f\omega^n\wedge\Phi
\end{eqnarray}

where the unknown function $\phi$ belongs to the class of $C^\infty$
functions such that $\omega-\sqrt{-1}\6\bar\6 \phi$ lies in the
interior of the cone of $\Phi$-positive forms. (i) Then on this
class of functions the Monge-Amp\`ere equation (\ref{E:ma1}) is
elliptic, and its solution is unique up to a constant. (ii) If
moreover the forms $\Phi$ and  $\omega\wedge\Phi$ are closed, then a
necessary condition of the solvability of (\ref{E:ma1}) is
$$\int_X(e^f-1)\omega^n\wedge\Phi=0.$$
\end{theorem}

{\bf Proof:} Part (i) follows immediately from \ref{L:4} and
\ref{L:7}. Let us prove part (ii). It is enough to show that, for
any $j$, one has
$$\int_X (\6\bar\6\phi)^j\wedge\omega^{n-j}\wedge\Phi=0.$$
This equality will follow from Stokes' formula if we prove that
$d(\omega^j\wedge\Phi)=0$ for any $j$. But
$$d(\omega^j\wedge\Phi)=j\omega^{j-1}\wedge
d\omega\wedge\Phi=j\omega^{j-1}\wedge d(\omega\wedge\Phi)=0.$$
Theorem is proven.
\endproof

\hfill

\begin{lemma}
Let $(M^{4n},I,J,K)$ be a hypercomplex manifold. Let $\Theta\in
\Lambda^{2n,0}_I(M)$ be q-real, q-positive, non-vanishing
holomorphic form. Let $\Phi=V(1)\in \Lambda^{n,n}_{I,+}(M)$ be as in
(\ref{D:Phi}). Let $\omega\in \Lambda^{1,1}_{I,+}(M)$ be a positive
form (in the complex sense). Then $\omega$ belongs to the interior
of the cone of $\Phi$-positive $(1,1)$-forms.
\end{lemma}

{\bf Proof:} This follows from \ref{L:2}, because $\Phi$ is weakly
positive by \ref{_P:Phi-character_Proposition_} and $\omega$ is
strongly positive by assumption. \endproof

\hfill

\begin{lemma}\label{L:interior}
Assume that $\phi$ satisfies the quaternionic Monge-Amp\`ere
equation $$(\Omega+\6\6_J\phi)^n=e^f\Omega^n$$ on a compact manifold
$M$. Then the form $\Omega+\6\6_J\phi$ belongs to the interior of
the cone of (strongly=weakly) q-positive $(2,0)$-forms. Hence
$\Omega+\6\6_J\phi$ is an HKT-form.
\end{lemma}

{\bf Proof:} Let $x\in M$ be a point where $\phi$ achieves its
minimum. One has $\6\6_J\phi(x)\geq 0$. Hence
$(\Omega+\6\6_J\phi(x))\geq 0$. But the top power
$(\Omega+\6\6_J\phi)^n$ is nowhere vanishing and continuous. Hence,
everywhere $\Omega+\6\6_J\phi$ belongs to the interior of the cone
of q-positive elements. \endproof

\hfill

\begin{corollary}\label{C:ellip-unique}
Let $(M^{4n},I,J,K)$ be a compact hypercomplex manifold. Let
$\Omega\in \Lambda^{2,0}_I(M)$ be an HKT-form. Let us assume
moreover that $M$ admits a non-vanishing holomorphic q-positive form
$\Theta\in \Lambda^{2n,0}_I(M)$. Fix a real-valued smooth function
$f$. Consider the quaternionic Monge-Amp\`ere equation
$$(\Omega+\6\6_J\phi)^n=e^f\Omega^n$$
on the class of $C^\infty$-smooth functions $\phi$. Then the
quaternionic Monge-Amp\`ere equation is elliptic and the solution is
unique up to a constant. Moreover a necessary condition for
solvability of this equation is
$$\int_M(e^f-1)\Omega^n\wedge\bar\Theta=0.$$
\end{corollary}

{\bf Proof:} The proof follows immediately from \ref{C:8},
\ref{L:interior} and the properties of the form $\Phi$ given in
\ref{_P:Phi-character_Proposition_}. \endproof

\hfill

\begin{remark}\label{_Holonomy_Remark_}
As we have already mentioned in the introduction, it was shown in
\cite{_Verbitsky:Canonical_HKT_} that if $M$ is a compact
HKT-manifold admitting a holomorphic (with respect to $I$)
$(2n,0)$-form $\Theta$ then the holonomy of the Obata connection is
contained in the group $SL_n(\HH)$ (instead of $GL_n(\HH)$).
Conversely, if the holonomy of the Obata connection is contained in
$SL_n(\HH)$ then there exists a form $\Theta$ as above which
moreover can be chosen to be q-positive (in sense of Section
\ref{Ss:R_V} below).
\end{remark}

\def\Om{\Omega}
\def\pfi{\pt\pt_J\phi}

\section[Zero order estimates for the quaternionic Monge-\-Amp\`ere
equation.]{Zero-order estimates for the quaternionic\\
Monge-\-Amp\`ere equation.}\label{S:zero-order}

In this section we will make the following assumption on an
HKT-manifold $M^{4n}$. We assume that $M$ is compact, connected, and
there exists a non-vanishing q-positive holomorphic section
$\Theta\in \Lambda^{2n,0}_I(M)$. The main result of this section is
\ref{ze-6}. Recall that we study the quaternionic Monge-Amp\`ere
equation
\begin{eqnarray}\label{ma10}
(\Omega_0+\6\6_J\phi)^n=e^f\Omega^n_0,
\end{eqnarray}
where $\phi$ is a real valued $C^\infty$-smooth function. By
\ref{L:interior}, $\Omega_0+\6\6_J\phi$ is an HKT-form.

Let us formulate a conjecture which is a quaternionic version of the
Calabi conjecture.

\hfill

\conjecture Let $(M^{4n},I,J,K)$ be a compact hypercomplex manifold
with an HKT-form $\Omega_0\in \Lambda^{2,0}_I(M)$. Assume, in
addition, that there exists a non-vanishing holomorphic q-positive
form $\Theta\in \Lambda^{2n,0}_I(M)$.\footnote{This is equivalent to
$\Hol(M)\subset SL(n, {\Bbb H})$, see \ref{_Holonomy_Remark_}} Then
the Monge-Amp\`ere equation (\ref{ma10}) has a $C^\infty$-solution
provided the following necessary condition is satisfied:
$$\int_M(e^f-1)\Ome_0^n\wedge\bar\Theta=0.$$

\hfill

Let $\phi\in C^2(M,\RR)$ be a solution of the Monge-Amp\`ere
equation
$$(\Omega_0+\pt\pt_J\phi)^n=e^f\Omega_0^n$$ satisfying the
normalization condition
\begin{equation}\label{_normali_CY_Equation_}
\int_M\phi\cdot \Omega^n_0\wedge\bar\Om^n_0=0.
\end{equation}


For brevity, we will denote $\Omega:=\Omega_0+\pt\pt_J\phi$. Let us
normalize the form $\Om_0$ such that $vol_{g_0}(M)=1$ where $g_0$ is
the HKT-metric corresponding to $\Omega_0$. The next lemma is
essentially linear algebraic.

\begin{lemma}\label{ze-2}
For any smooth function $\psi$ one has pointwise
$$|\nabla \psi|^2_{g_0}=4n \cdot \frac{\pt\psi \wedge \pt_J\psi \wedge
\Om_0^{n-1}}{\Om_0^n},$$ where $|\cdot|_{g_0}$ denotes the norm on
$TM$ with respect to $g_0$.
\end{lemma}

{\bf Proof:} The proof is elementary and is left to a reader.
\endproof

\hfill

\begin{proposition}\label{ze-1}
Let $p>1$. Then the solution $\phi$ satisfies the following estimate
$$||\nabla|\phi|^{p/2}||_{L^2}^2\leq
\frac{1}{16n}\cdot\frac{p^2}{(p-1)} \int_M(1-e^f)\phi
|\phi|^{p-2}\Om_0^n\wedge\bar\Theta.$$
\end{proposition}

{\bf Proof:} We have
\begin{eqnarray*}
\int_M(1-e^f)\phi |\phi|^{p-2}\Om_0^n\wedge\bar\Theta =\int_M \phi
|\phi|^{p-2}(\Om_0^n-\Om^n)\wedge\bar\Theta=\\
-\int_M \phi |\phi|^{p-2}\pfi\wedge(\sum_{l=0}^{n-1}\Om_0^l\wedge
\Om^{n-1-l}) \wedge\bar\Theta=\\\int_M\pt(\phi
|\phi|^{p-2})\wedge\pt_J\phi\wedge(\sum_{l=0}^{n-1}\Om_0^l\wedge
\Om^{n-1-l})\wedge\bar\Theta=\\
(p-1)\int_M |\phi|^{p-2}\pt\phi \wedge
\pt_J\phi\wedge(\sum_{l=0}^{n-1}\Om_0^l\wedge
\Om^{n-1-l})\wedge\bar\Theta.
\end{eqnarray*}
Since $\Om_0,\Om,$ and $\Theta$ are positive the last expression is
at least
$$(p-1)\int_M |\phi|^{p-2}\pt\phi \wedge
\pt_J\phi\wedge\Om_0^{n-1}\wedge\bar\Theta.$$ But
$$\pt|\phi|^{p/2}=\frac{p}{2}|\phi|^{p/2-1}\pt \phi,
\pt_J|\phi|^{p/2}=\frac{p}{2}|\phi|^{p/2-1}\pt_J \phi.$$ Thus we get
\begin{eqnarray*}
\int_M(1-e^f)\phi |\phi|^{p-2}\Om_0^n\wedge\bar\Theta \geq
(p-1)\int_M
|\phi|^{p-2}\pt\phi \wedge \pt_J\phi\wedge\Om_0^{n-1}\wedge\bar\Theta\geq\\
(p-1)\frac{4}{p^2}\int_M
\pt|\phi|^{p/2}\wedge\pt_J|\phi|^{p/2}\wedge\Om_0^{n-1}\wedge\bar\Theta\overset{\ref{ze-2}}{=}
16n\cdot\frac{p-1}{p^2}|\nabla|\phi|^{\frac{p}{2}}|^2_{g_0}.
\end{eqnarray*}
This implies \ref{ze-1}. \endproof

\hfill

In this section we use the notation $\kappa:=\frac{2n}{2n-1}$. Let
us denote by $L^2_1(M)$ the Sobolev space of functions on $M$ such
that all partial derivatives up to order 1 are square integrable.

\hfill

\begin{lemma}\label{ze-3}
There exists a constant $C_1$ depending on $M$ and $\Omega_0$ only,
such that, for any function $\psi\in L_1^2(M)$,
$$||\psi||^2_{L^{2\kappa}}\leq C_1(||\nabla
\psi||_{L^2}^2+||\psi||^2_{L^2}).$$ Moreover, if  $\psi$ satisfies
$\int_M\psi\cdot \Om_0^n\wedge\bar\Om^n_0=0$, one has
$$||\psi||^2_{L^{2\kappa}}\leq C_1||\nabla \psi||_{L^2}^2.$$
\end{lemma}

{\bf Proof:} By the Sobolev imbedding theorem there exists a
constant $C'$ such that, for any function $\psi\in L_1^2(M)$, one
has
$$||\psi||^2_{L^{2\kappa}}\leq C'(||\nabla
\psi||^2_{L^2}+||\psi||^2_{L^2}).$$  If the function $\psi$
satisfies $\int_M\psi\cdot \Om_0^n\wedge\bar\Om^n_0=0$ then one has
$||\psi||_{L^2}^2\leq \tilde C ||\nabla \psi||_{L^2}^2$ since the
second eigenvalue of the Laplacian on $M$ is strictly positive. Thus
\ref{ze-3} is proven.
\endproof

\hfill

\begin{lemma}\label{ze-4}
There exists a constant $C_2$ depending on $M,g_0,$ and
$||f||_{C^0}$ only such that if $p\in[2,2\kappa]$ then
$||\phi||_{L^p}\leq C_2$.
\end{lemma}

{\bf Proof:} Let us put $p=2$ in \ref{ze-1}. We get
\begin{align*}
 ||\nabla \phi||^2_{L^2}\leq & 4 \const
 \exp({||f||_{C^0}})||\phi||_{L^1}\\ \leq &4 \const\cdot
 vol_{g_0}(M)^{1/2}
 \exp({||f||_{C^0}})||\phi||_{L^2}
\end{align*}
where the second
 inequality follows from the H\"older inequality. Since
\[ \int_M\phi\cdot\Om_0^n\wedge\bar\Om^n_0=0,\] we have
 $$||\phi||_{L^2}\leq C||\nabla \phi||_{L^2}.$$
Hence $||\nabla \phi||_{L^2}\leq C\cdot 4 const\cdot
vol_{g_0}(M)^{1/2} \exp({||f||_{C^0}})$. Therefore by \ref{ze-3}
there exists a constant $C_2'$ depending on $M,g_0,$ and
$||f||_{C^0}$ only such that $$||\nabla\phi||_{L^{2\kappa}} \leq
C_2'.$$ Hence by the H\"older inequality $||\nabla\phi||_{L^{p}}\leq
C_2''$ for $p\in[2,2\kappa]$. \endproof

\begin{proposition}\label{ze-5}
There exist constants $Q_1,C_3$ depending on $M,$ $g_0,$
$||f||_{C^0}$ only such that for any $p\geq 2$
\[ ||\phi||_{L^p}\leq Q_1(C_3p)^{-\frac{2n}{p}}.\]
\end{proposition}

{\bf Proof:} Define $C_3=C_1(2 \cdot const \cdot  e^{||f||_{C^0}}
+1)\cdot \kappa^{(2n-1)}$ where $const$ is from \ref{ze-1}. Choose
$Q_1$ so that $Q_1>C_2(C_3p)^{\frac{2n}{p}}$ for $2\leq p\leq
2\kappa$ and $Q_1>(C_3p)^{\frac{2n}{p}}$ for $2\leq p<\infty$. We
will prove the result by induction on $p$. By \ref{ze-4}, if $2\leq
p\leq 2\kappa$ then $||\phi||_{L^p}\leq C_2\leq
Q_1(C_3p)^{-\frac{2n}{p}}$ . For the inductive step, suppose that
$$||\phi||_{L^p}\leq Q_1(C_3p)^{-\frac{2n}{p}} \mbox{ for } 2\leq p\leq
k,\, \mbox{ where } k\geq 2\kappa \mbox{ is a real number}.$$ We
will show that, for
$$||\phi||_{L^q}\leq Q_1(C_3q)^{-\frac{2n}{q}} \mbox{ for } 2\leq
q\leq\kappa k,$$ and therefore by induction \ref{ze-5} will be
proved. Let $p\in[2,k]$. By \ref{ze-1} we get
\begin{eqnarray}\label{ine1}
||\nabla|\phi|^{p/2}||_{L^2}^2\leq
const\frac{p^2}{(p-1)}e^{||f||_{C^0}}||\phi||^{p-1}_{L^{p-1}}.
\end{eqnarray}

Applying \ref{ze-3} to $\psi=|\phi|^{p/2}$ we get
\begin{eqnarray}\label{ine2}
||\phi||^p_{L^{\kappa p}}\leq
C_1(||\nabla|\phi|^{p/2}||^2_{L^2}+||\phi||^p_{L^p}).
\end{eqnarray}

Combining (\ref{ine2}) and (\ref{ine1}) we obtain
$$||\phi||^p_{L^{\kappa p}}\leq
C_1(2p \cdot const \cdot
e^{||f||_{C^0}}||\phi||^{p-1}_{L^{p-1}}+||\phi||^p_{L^p}).$$ Let
$q=\kappa p$. Since $2\leq p\leq k$ we have $||\phi||_{L^p}\leq
Q_1(C_3p)^{-\frac{2n}{p}}$. Since $||\phi||_{L^{p-1}}\leq
||\phi||_{L^p}$ we get
\begin{eqnarray*}
||\phi||^p_{L^q}\leq  C_1\left(2p \cdot const \cdot
e^{||f||_{C^0}}||\phi||^{p-1}_{L^{p}}+
\left(Q_1(C_3p)^{-\frac{2n}{p}}\right)^p\right)\leq\\
C_1\left(2p \cdot const \cdot
e^{||f||_{C^0}}(Q_1(C_3p)^{-\frac{2n}{p}})^{p-1}+\left(Q_1(C_3p)^{-\frac{2n}{p}}\right)^p\right).
\end{eqnarray*}

But $Q_1(C_3p)^{-\frac{2n}{p}}\geq 1$. Hence
$$||\phi||^p_{L^q}\leq C_1Q_1^p(C_3p)^{-2n}(2 \cdot const \cdot p
e^{||f||_{C^0}} +1).$$ It remains to show that the last expression
is at most $Q_1^p(C_3q)^{-\frac{2n}{q} p}.$ It is enough to check
that $$C_1(C_3p)^{-2n}(2p \cdot const \cdot e^{||f||_{C^0}} +1)\leq
(C_3q)^{-\frac{2n}{\kappa}}=(C_3\kappa p)^{-(2n-1)}.
$$
The left-hand side is at most $C_1(C_3p)^{-2n}\cdot p(2 \cdot const
\cdot  e^{||f||_{C^0}} +1)$. Hence it is enough to check that
$$C_1C_3^{-2n}(2 \cdot const \cdot  e^{||f||_{C^0}} +1)\leq
(C_3\kappa )^{-(2n-1)}.$$ Namely
$$C_1(2 \cdot const \cdot  e^{||f||_{C^0}} +1)\leq C_3\cdot
\kappa^{-(2n-1)}.$$ But this holds by the definition of $C_3$.
\endproof

\hfill

The following corollary is the main result of this section.

\hfill

\begin{corollary}\label{ze-6}
The exists a constant $C_4$ depending on $M,g_0,||f||_{C_0}$ only,
such that $$||\phi||_{C^0}\leq C_4,$$ for any solution of
quaternionic Calabi-Yau equation \eqref{ma10} which satisfies the
normalization condition \eqref{_normali_CY_Equation_}.
\end{corollary}

{\bf Proof:} We have $$||\phi||_{C^0}=\lim_{p\to
\infty}||\phi||_{L^p}\leq Q_1$$ where the last inequality follows
from \ref{ze-5}. \endproof



\hfill

\hfill

\hfill

\noindent {\sc Semyon Alesker\\
{ \normalsize Department of Mathematics, Tel Aviv University, Ramat
Aviv}
 \\  { \normalsize 69978 Tel Aviv,
Israel }
\\ \tt semyon@post.tau.ac.il\\

\hfill

\noindent {\sc Misha Verbitsky\\
{\sc  Institute of Theoretical and
Experimental Physics \\
B. Cheremushkinskaya, 25, Moscow, 117259, Russia }\\
\tt verbit@maths.gla.ac.uk, \ \  verbit@mccme.ru
}
}
\end{document}